\documentclass[11pt]{amsart}
\usepackage[utf8x]{inputenc}
\usepackage[a4paper]{geometry}
\usepackage{amsthm,amstext,amsmath,amscd,amssymb,latexsym}
\usepackage{hyperref}
                 \hypersetup{ pdfborder={0 0 0}, 
                              colorlinks=true, 
                              linktoc=page,
                              pdfauthor={M. C. Brambilla and G. Staglian{\`o}}, 
                              pdftitle={On the algebraic boundaries among typical ranks}} 
\usepackage{color}
\usepackage[dvipsnames]{xcolor}
\usepackage{fancyvrb}   
 
\definecolor{airforceblue}{rgb}{0.36, 0.54, 0.66}
\definecolor{bleudefrance}{rgb}{0.19, 0.55, 0.91}
\definecolor{darkorchid}{rgb}{0.6, 0.2, 0.8}
\definecolor{darkorange}{rgb}{1.0, 0.55, 0.0}
\definecolor{darkspringgreen}{rgb}{0.09, 0.45, 0.27}

\newcommand{\CC}{{\mathbb C}}
\newcommand{\PP}{{\mathbb P}}

\newcommand{\RR}{{\mathbb R}}

\newcommand{\Sym}{{\textrm{Sym}}}
\newcommand{\Sing}{{\rm{Sing}}}
\newcommand{\rk}{{\rm{rk}}}

\DeclareMathOperator{\CH}{CH}
\DeclareMathOperator{\Gr}{\mathbb{G}}

\newtheorem{theorem}{Theorem}[section]
\newtheorem{lemma}[theorem]{Lemma}
\newtheorem{proposition}[theorem]{Proposition}
\newtheorem{corollary}[theorem]{Corollary}

\theoremstyle{definition}

\newtheorem{example}[theorem]{Example}

\theoremstyle{remark}
\newtheorem{remark}[theorem]{Remark}

\numberwithin{equation}{section}

\title[On the algebraic boundaries among typical ranks]{On the algebraic boundaries among typical ranks for real binary forms}

\author[M. C. Brambilla]{Maria Chiara Brambilla}
\email{{\tt brambilla@dipmat.univpm.it}}
\address{Universit\`a Politecnica delle Marche, 
  via Brecce Bianche, I-60131 Ancona, Italy}

\author[G. Staglian\`o]{Giovanni Staglian\`o}
\email{{\tt giovannistagliano@gmail.com}}
\address{Universit\`a Politecnica delle Marche, 
   via Brecce Bianche, I-60131 Ancona, Italy}

\thanks{The first author is partially supported by MIUR and INDAM}



\begin{document}

\begin{abstract}
We describe the algebraic boundaries of the regions of real binary forms with fixed typical rank and of degree at most eight, 
showing that they are dual varieties of suitable coincident root loci.
\end{abstract}

\maketitle

\section{Introduction}

The study of 
symmetric tensors, of their rank, decomposition and identifiability is a classical problem, 
which received great attention recently in both pure and applied mathematics; 
see e.g.\ \cite{Landsberg} and references therein, 
see also \cite{BERNARDI2014205,BBO,Michalek2017,Carlini2017,MASSARENTI2018227,RealComplex}.

Symmetric tensors can be interpreted as homogeneous polynomials, also called forms.
The 
rank of a degree $d$ form $f$ is 
the minimum integer $r$ such that there exists a decomposition 
$f=\sum_{i=1}^rc_i(l_i)^d$, where $l_i$ are linear forms and 
$c_i$ are scalars.

In this paper we focus on the case of binary forms 
over the field of real numbers $\mathbb{R}$.
In this case it is known that the (real) rank of a general form satisfies the inequalities 
$\frac {d+1}2\le r\le d$. Moreover all the ranks in this range are typical, 
that is, they occur in open 
subsets (with respect to the Euclidean topology) of the real vector space of degree $d$ forms; see \cite{Blekherman}.

A natural problem is to understand the geometry of the sets $\Omega_{d,r}$ of forms of degree $d$ and rank $r$.
In particular we would like to describe the boundaries among the various sets of forms of given typical rank;
more precisely,
we are interested in understanding the algebraic boundaries,
i.e., 
the Zariski closures of the topological boundaries 
(see Section~\ref{secBinForms} for the precise definitions).

The easiest case is the maximal one, that is when the rank is equal to the degree. Indeed it is proved in 
\cite{CausaRe, ComonOttaviani} that a binary form of degree $d$ 
with distinct roots has rank $d$ if and only if all its roots are real. 
Hence its algebraic
boundary is the discriminant hypersurface of forms with two coincident roots.

The geometric description of the sets $\Omega_{d,r}$ becomes much more intricate for $r<d$. Indeed, 
although the rank of a form is always greater than or equal to the number of its real distinct roots, 
in general the number of real distinct roots is not invariant in the region 
$\Omega_{d,r}$.

In \cite{LeeSturmfels} the authors study the boundary of the set of forms of rank 
$\lceil\frac {d+1}2\rceil$, which is the minimal typical rank. They
prove that the components of the boundary are dual varieties of suitable coincident root loci.

We tackle the problem of describing all the intermediate boundaries in general, 
as proposed by Lee and Sturmfels in \cite[Remark 4.5]{LeeSturmfels}.
Our approach provides a unified description of all the boundaries
in terms of dual varieties of coincident root loci.
We recall that the cases of degree $d\le5$ have been described in 
\cite{ComonOttaviani}, while the case $d=6$ follows by \cite{LeeSturmfels,CausaRe} 
(see Proposition \ref{propLeq6}  for more details).
In this paper we focus on the cases  $d=7$ and $d=8$, and we postpone a general description to future work. 

The paper is organized as follows. Sections \ref{secCRL} and \ref{secBinForms} 
are devoted to preliminary results; in particular,
in Proposition~\ref{propLeq6} we recall the known results concerning algebraic boundaries for real binary forms 
of degree less than or equal to $6$.
Section \ref{secCaso7} and \ref{secCaso8} contain our main results,
which are Theorem \ref{teorema} and Theorem \ref{teorema2}, describing 
the  algebraic
boundaries for real binary forms of degrees respectively  $7$ and $8$. They turn out to be dual varieties
of suitable coincident root loci.
Finally in Section \ref{secComp} we explain some of the computational methods of which 
we take advantage in our study.

\subsection*{Acknowledgements} 
We thank the participants to the seminar ``Algebraic Geometry and Tensors'' where this work has begun. 
In particular we are grateful to Giorgio Ottaviani for many useful discussions.

\section{Coincident root loci}\label{secCRL}
We recall here some known results on coincident root loci, referring to 
\cite{Weyman-1, Katz, Chipalkatti-1, Chipalkatti-2, Kurmann}
for details.

We regard a degree $d$ binary form $f=\sum_{i=0}^d \binom{d}{i} a_i x^{d-i} y^i$ 
over the complex field $\CC$
as a point of the projective space $\PP(\CC[x,y]_d)$, where
$\CC[x,y]_d=\Sym^d(\CC^2)$. This space is identified with 
$\PP^d$ 
using homogeneous coordinates $(a_0,\ldots,a_d)$.

A partition $\lambda=(\lambda_1,\ldots,\lambda_n)$ 
of $d$ is a list of integers $\lambda_1\ge{\cdots\ge}\lambda_n\ge1$ such that $\sum_{i=1}^n \lambda_i=d$.
Given a partition $\lambda$,
the {\it coincident root locus}
$\Delta_\lambda\subset \PP^d$ is the
set of binary forms $f$ of degree $d$ 
{which admit a factorization $f=\prod_{i=1}^n \ell_i^{\lambda_i}$}
for some linear forms $\ell_1,\ldots,\ell_n\in \CC[x,y]_1$.
A partition $\lambda$ can be also represented by the list of integers $m_1,\ldots, m_k$ 
defined as $m_j=|\{i:\lambda_i=j\}|$, and clearly $\sum_{j=1}^k jm_j=d$.
Then 
the coincident root locus $\Delta_\lambda$ is given by the binary forms of degree $d$ 
which have $m_j$ roots of multiplicity at least $j$.  
It is classically known (see \cite{Hilbert}) that $\Delta_\lambda\subset\PP^d$ is a variety of dimension $n$ and degree
\begin{equation}
\deg(\Delta_\lambda)=\frac{n!}{m_1!m_2!\cdots m_k!}\lambda_1\lambda_2\cdots\lambda_n  .
\end{equation}
If $\lambda=(2,1^{d-2})$, the corresponding coincident root locus $\Delta_\lambda=\Delta$ 
is the classical discriminant  hypersurface.
In the opposite case, if $\lambda=(d)$ then $\Delta_\lambda$ is the rational normal curve $C_d\subset\PP^d$. 
When $\lambda=(a,1^{d-a})$, the partition is called {\it hook}, and the associated
coincident root locus $\Delta_\lambda$ 
 represents the tangential developable of $\Delta_{(a+1,1^{d-a-1})}$.

\subsection{Singularities of \texorpdfstring{$\Delta_\lambda$}{Dl}}
The singular loci of coincident root loci have been studied by Chipalkatti \cite{Chipalkatti-1} and Kurmann \cite{Kurmann}. 

Given a partition $\lambda=(\lambda_1,\ldots,\lambda_n)$, the singular locus $\Sing(\Delta_\lambda)$  
is given by the union of $\Delta_{\mu}$ for some suitable coarsenings $\mu$ of $\lambda$. 
See either \cite[Definition 5.2]{Chipalkatti-1}, 
or \cite[Proposition 2.1]{Kurmann} for the precise description.
In particular $\Delta_\lambda$ is smooth if and only if $\lambda_1=\cdots=\lambda_n$.
Otherwise the singular locus is of {(not necessarily pure)}
codimension $1$. 

\begin{example}\label{esempioSING}
For future use, we now
 compute the \emph{ iterate singular locus} of $\Delta_\lambda$, for $\lambda=(2,1,1,1)$ and $\lambda=(2,1,1,1,1)$.
 \begin{gather*}
 \Sing(\Delta_{(2,1,1,1)}) = \Delta_{(3,1,1)}\cup \Delta_{(2,2,1)} ; \\
 \Sing(\Delta_{(3,1,1)}) = \Delta_{(4,1)}, \quad 
 \Sing(\Delta_{(2,2,1)}) = 
 \Delta_{(3,2)} ; \\ 
 \Sing(\Delta_{(4,1)}) = \Sing(\Delta_{(3,2)}) = \Delta_{(5)} .
 \end{gather*}
 \begin{gather*}
 \Sing(\Delta_{(2,1,1,1,1)}) = \Delta_{(3,1,1,1)}\cup \Delta_{(2,2,1,1)} ; \\
 \Sing(\Delta_{(3,1,1,1)}) = \Delta_{(4,1,1)}\cup \Delta_{(3,3)}, \quad 
 \Sing(\Delta_{(2,2,1,1)}) = \Delta_{(3,2,1)} \cup \Delta_{(2,2,2)} ;\\
 \Sing(\Delta_{(4,1,1)}) = \Delta_{(5,1)}, \quad 
 \Sing(\Delta_{(3,2,1)}) =  \Delta_{(3,3)}\cup \Delta_{(4,2)}\cup \Delta_{(5,1)}  ; \\
 \Sing(\Delta_{(4,2)}) = \Sing(\Delta_{(5,1)}) = \Delta_{(6)} .
 \end{gather*}
\end{example}

\subsection{Duality} \label{duality-sec}

Consider  the dual ring of differential operators
 $\CC[\partial_x,\partial_y]=\CC[u,v]$, which acts on $\CC[x,y]$ 
 with
the usual rules of differentiations
 and gives the pairing with respect to the degrees,
$$\CC[x,y]_d\otimes \CC[u,v]_k\to \CC[x,y]_{d-k} .$$
The conormal variety of a coincident root locus $\Delta_\lambda$ is
 the Zariski closure of the set 
$$\{(f,g):f\mbox{ is a smooth point of $\Delta_\lambda$}, g\perp T_f\Delta_\lambda\}\subset \PP(\CC[x,y]_d)\times\PP(\CC[u,v]_d)  ,  $$
where
$T_f\Delta_\lambda$
denotes the tangent space to $\Delta_\lambda$ at a point $f$.
The dual variety  $(\Delta_\lambda)^\vee$ of $\Delta_\lambda$ is 
 the projection onto $\PP(\CC[u,v]_d)$ of the conormal variety of $\Delta_\lambda$. 
The biduality theorem (see \cite{GKZ}) implies that $(\Delta_\lambda^\vee)^\vee=\Delta_\lambda$. 

Lee and Sturmfels study duality for binary forms in \cite{LeeSturmfels}. 
We recall here some results which we will use in the sequel.

\begin{proposition}[\cite{LeeSturmfels}]
\label{proposition-LS}
Given $\lambda=(\lambda_1,\ldots,\lambda_n)$ and $\Delta_\lambda\subset\PP(\CC[u,v]_d)$,
the points of the dual variety $\Delta_\lambda^\vee\subset \PP(\CC[x,y]_d)$
are given by the binary forms $f(x,y)$ that are annihilated by some order $d-n$ operator of the form 
$\Pi_{i=1}^n \ell_i^{\lambda_i-1}(\partial_x,\partial_y)$ where $\ell_i\in \CC[u,v]_1$.
\end{proposition}

\begin{proposition}[\cite{LeeSturmfels}]\label{propLSJoin}
\label{dual-hook}
Given  $\lambda=(\lambda_1,\ldots,\lambda_n)$ and $\Delta_\lambda\subset\PP(\CC[u,v]_d)$, 
the dual variety $\Delta_\lambda^\vee\subset \PP(\CC[x,y]_d)$ has codimension 
$m_1 + 1$, 
and
it is given by the join of the $(n-m_1)$ coincident root loci 
$\Delta_{(d-\lambda_i+2,1^{\lambda_i-2})}$ for $1\le i\le n$ with $\lambda_i\ge2$.
\end{proposition}
If $\lambda_i\ge2$ for all $i$, then 
$\Delta_\lambda^\vee$ is a hypersurface 
of degree (see \cite{Oeding})
\begin{equation}\label{luke}
\frac{(n+1)!}{m_2!\cdots m_k!}(\lambda_1-1)(\lambda_2-1)\cdots(\lambda_n-1).
\end{equation}

\subsection{Chow forms and higher associated varieties} 

Let $\Gr(h,m)$ denote the Grassmannian of projective subspaces of dimension $h$ in $\PP^m$.
Let $X\subset \PP^m$ be a projective variety  of dimension $k$.

The {\it $i$-th higher associated variety} $\CH_i(X)$ of $X$ is defined as 
the closure of the set of all $(m-k-1+i)$-dimensional subspaces 
$L\subset \PP^m$ such that $L\cap X\neq \emptyset$ and $\dim(L\cap T_xX)\ge i$ 
for some smooth point $x\in L\cap X$ (where $T_xX$ denotes
 the embedded tangent space to $X$ at $x$), see \cite{GKZ} for details.

For $i=0$, the associated variety $\CH_0(X)\subset\Gr(m-k-1,m)$ is the Chow hypersurface, 
while for $i=k$, we have that $\CH_k(X)\subset\Gr(m-1,m)$ corresponds to the dual variety $X^\vee$ 
via the Grassmannian duality $\Gr(m-1,\PP^m)\simeq\Gr(0,(\PP^m)^\vee)$.
If $i=1$ and $\deg(X)\ge2$, the associated variety 
$\CH_1(X)$ is the Hurwitz hypersurface,  see \cite{Sturmfels-Hurw}.

The variety $\CH_i(X)$ is a hypersurface
if and only if $i\le \dim(X)-(m-1- \dim(X^\vee))$, 
see \cite{kohn}.
In particular, if $X=\Delta_{(\lambda_1,\ldots,\lambda_n)}$ 
is a coincident root locus, the higher associated variety $CH_i(X)$ is a hypersurface if
and only if $i\le|\{j:\lambda_j\ge2\}|$.
 
\section{Real rank of binary forms}\label{secBinForms}

\subsection{Typical ranks for binary forms}
Given a binary form $f$ of degree $d$ with complex (or real) coefficients, 
its {\it complex rank} is
the minimum integer $r$ such that $f$ admits a decomposition $f=\sum_{i=1}^r (\ell_i)^d$ 
where $\ell_i$ are linear forms with complex coefficients.
The {\it generic complex rank} for binary forms of degree $d$ (that is the rank of a general binary form of degree $d$)
is $\lceil\frac {d+1}2\rceil$.
Sylvester Theorem says that a general binary form admits a unique {minimal} decomposition if the degree is 
odd, infinitely many (parametrized by a line) if the degree is even.

Consider now the polynomial ring $R=\RR[x,y]$ of real binary forms.
Given $f\in R_d$, the {\it real rank} of  $f$ (denoted by $\rk(f)$)
is the minimum integer $r$ such that $f$ admits a decomposition $f=\sum_{i=1}^r c_i (l_i)^d$ 
where $l_i\in R_1$ and
$c_i\in\mathbb{R}$;
we can impose $c_i\in\{1,-1\}$ if $d$ is even, and $c_i=1$ if $d$ is odd.

In the real field the notion of generic rank is replaced by the notion of {\it typical ranks}. 
A rank is called typical for binary forms of degree $d$ if it occurs in an open subset of $R_d$, 
with respect to the Euclidean topology.

Define
$\Omega_{d,r}=\{f\in R_d: \rk(f)=r\}$,
and denote by 
$\mathcal{R}_{d,r} 
$ the interior of $\Omega_{d,r}$.
Then
$\mathcal{R}_{d,r}$ is a semi-algebraic set in the real vector space $R_d$, and a 
rank is typical exactly when $\mathcal{R}_{d,r}$ is not empty. 
From the main result of \cite{Blekherman}, 
a rank $r$ is typical if and only if $\frac {d+1}2\le r \le d$.
Thus, from now on we assume that  $\frac {d+1}2 \le r \le d$.
We define the {\it topological boundary} $\partial(\mathcal{R}_{d,r})$ as the set-theoretic difference 
of the closure of $\mathcal{R}_{d,r}$
and 
the interior of the closure of $\mathcal{R}_{d,r}$. 
It is a semi-algebraic subset of $R_d$ of pure codimension one. 
We define the \emph{real rank boundary}
$\partial_{\textrm{alg}}(\mathcal{R}_{d,r})$
as the Zariski closure of the topological boundary  $\partial(\mathcal{R}_{d,r})$ (see also \cite[Section~4]{LeeSturmfels}).
The real rank boundaries $\partial_{\textrm{alg}}(\mathcal{R}_{d,r})$ are hypersurfaces 
of the real 
space $R_d$, that we consider as 
hypersurfaces of the complex projective space $\PP(\CC[x,y]_d)=\PP_\CC^d$.
Let us remark that these hypersurfaces are 
invariant under the natural action of $\mathrm{SL}_2$ on $\PP^d$.

Real rank boundaries have been studied only in the two extreme cases, that is for
 maximum rank $d$ and minimum rank $\overline{r}=\lceil\frac {d+1}2\rceil$.
In the first case $\partial_{\textrm{alg}}(\mathcal{R}_{d,d})$
is the discriminant hypersurface $\Delta_{(2,1^{d-2})}$ (see 
\cite[Proposition~3.1]{ComonOttaviani} and \cite[Corollary~1]{CausaRe});
in the second case the real rank boundary 
$\partial_{\textrm{alg}}(\mathcal{R}_{d,\overline{r}})$
is described in \cite[Theorem 4.1]{LeeSturmfels}.
Hence, for $d\le 6$ we have  a complete description of all the real rank boundaries, 
that we recall in the following:
\begin{proposition}[\cite{ComonOttaviani,CausaRe,LeeSturmfels}]\label{propLeq6}
 The real rank boundaries for binary forms of degree $\leq6$ are the following hypersurfaces:
   \begin{align*}
    \partial_{\textrm{alg}}(\mathcal{R}_{3,2}) &= \partial_{\textrm{alg}}(\mathcal{R}_{3,3}) = (\Delta_{(3)})^\vee; 
    \end{align*}
    \begin{align*}
  \partial_{\textrm{alg}}(\mathcal{R}_{4,3}) &= \partial_{\textrm{alg}}(\mathcal{R}_{4,4}) = (\Delta_{(4)})^\vee;
  \end{align*}
  \begin{align*}
   \partial_{\textrm{alg}}(\mathcal{R}_{5,3}) &=(\Delta_{(3,2)})^\vee , \\
  \partial_{\textrm{alg}}(\mathcal{R}_{5,4}) &=(\Delta_{(3,2)})^\vee\cup (\Delta_{(5)})^\vee ,\\
  \partial_{\textrm{alg}}(\mathcal{R}_{5,5}) &=(\Delta_{(5)})^\vee ;
  \end{align*}
  \begin{align*}
   \partial_{\textrm{alg}}(\mathcal{R}_{6,4}) &= (\Delta_{(3,3)})^\vee\cup (\Delta_{(4,2)})^\vee , \\ 
   \partial_{\textrm{alg}}(\mathcal{R}_{6,5}) &=(\Delta_{(3,3)})^\vee\cup (\Delta_{(4,2)})^\vee \cup (\Delta_{(6)})^\vee  ,\\
  \partial_{\textrm{alg}}(\mathcal{R}_{6,6}) &=(\Delta_{(6)})^\vee . 
  \end{align*}
\end{proposition}
\begin{remark}
The hypersurfaces $(\Delta_{(3)})^\vee$, $(\Delta_{(4)})^\vee$, $(\Delta_{(5)})^\vee$, $(\Delta_{(6)})^\vee$
coincide with the discriminant hypersurfaces for binary forms of degrees $3,4,5,6$ and have degrees
$4,6,8,10$, respectively. For the other components, we have
\begin{itemize}
\item $(\Delta_{(3,2)})^\vee = \mathrm{Join}(\Delta_{(4,1)},\Delta_{(5)})$
is a hypersurface of  degree $12$ (this is the \emph{apple invariant $I_{12}$} considered in \cite{ComonOttaviani});
\item $(\Delta_{(3,3)})^\vee = \mathrm{Join}(\Delta_{(5,1)},\Delta_{(5,1)})$
is a hypersurface of  degree $12$;
\item $(\Delta_{(4,2)})^\vee = \mathrm{Join}(\Delta_{(4,1,1)},\Delta_{(6)})$
is a hypersurface of  degree $18$.
\end{itemize}
\end{remark}
\subsection{Apolarity}
We recall here classical techniques, going back to Sylvester.
Even if the results of this section are more general, we present them in the case of real numbers.
Let $R=\RR[x,y]$ be the polynomial ring of real binary forms 
and let  $D=\RR[\partial_x,\partial_y]=\RR[u,v]$ be the corresponding dual ring.
Given $l=ax+by \in R_1$,
the {\it apolar operator} is $l^\perp=-b\partial_x+a\partial_y\in D_1$.
Given a form $f$ in $R_d$, the {\it apolar ideal} $f^\perp\subset D$ is given by all the operators which annihilates $f$, that is:
$f^\perp=\{g(\partial_x,\partial_y)\in D:g\perp f \}$.
A basic tool is the following:
\begin{lemma}[Apolarity lemma]\label{apolarityLemma}
\label{apolarity}
  Assume $f\in R_d$ and let $l_i\in R_1$ be distinct linear forms for $1\le i\le r$. 
There are coefficients $c_i\in \RR$ such that $f=\sum_{i=1}^rc_i(l_i)^d$ if and only 
if the operator  $l_1^\perp\circ \cdots \circ l_r^\perp$ is in the apolar ideal $f^\perp$.
  \end{lemma}
We will say that a form of degree $d$ is {\it real-rooted} if it admits $d$ distinct real  roots.
From Lemma~\ref{apolarityLemma}, it follows that
a form $f$ has rank less than or equal to $r$ if and only if $(f^\perp)_r = f^\perp \cap D_r$ contains a 
real-rooted form. So the rank of $f$ is the smallest degree $r$ such that $(f^\perp)_r$ contains a 
real-rooted form.
The following result is an elementary
consequence of Lemma~\ref{apolarityLemma}. 
\begin{corollary}\label{rangoreale}
   Let $f$ be a real binary form, and let $r$ be an integer. 
 Then $\rk(F) < r$ if and only if $(f^\perp)_r\subset D_r$ 
   contains a special line whose generic member is a real-rooted form.
   Here, we say that a line $\langle g, g'\rangle \subset D_r$ 
   is \emph{special}
   if $\gcd(g,g')$ is a form of degree $r-1$.
\end{corollary}
The space of operators of degree $r$ contained in $f^\perp$ 
is the kernel of the linear map $A_{f}:D_r\to R_{d-r}$.
The 
{\it catalecticant (or Hankel) matrix} of $f$ is the matrix $A_{f}^{d,r}$ of size $(d-r+1)\times(r+1)$
that represents $A_f$ with respect to the standard basis.
We denote by $A^{d,r}$ the generic catalecticant matrix of size $(d-r+1)\times(r+1)$.

The following result is well-known (see e.g.\ \cite{IarrobinoKanev}):
\begin{proposition}
  Assume that $f\in R_d$ has rank greater than or equal to $2$.
  Then its apolar ideal $f^\perp$ is 
  generated by two real forms $g,g'$ such that $\deg g+\deg g'=d+2$ and
  $\mathrm{gcd}(g,g') = 1$.
  Conversely, any two such forms generate an ideal $f^\perp$ for some $f\in R$ with degree $\deg g+\deg g'-2$.  
\end{proposition}
We say that $f\in R_d$ is {\it generated in generic degrees} 
if  
  $(\deg g,\deg g') = (\lceil\frac {d+1}2\rceil,\lfloor\frac {d+3}2\rfloor)$.
The forms that are not generated in generic degrees form a subvariety of $R_d$.
More precisely, when the degree $d=2k$ is even, it is 
the hypersurface defined by the determinant of the intermediate $(k+1)\times (k+1)$
catalecticant matrix $A^{d,k}$; when the degree $d=2k+1$
 is odd, it is the subvariety of codimension $2$ defined by the maximal minors
 of the intermediate $(k+1)\times (k+2)$ catalecticant matrix $A^{d,k+1}$.

If $f$ is a binary form of degree $d$, with $\frac{d}{2}\le r\le d$ and having 
catalecticant matrix $A_{f}^{d,r}$ of maximal rank, then
$\dim(f^\perp)_r=2r-d$.
Thus, we can consider the \emph{apolar map}
$$\varPsi_{d,r}: \PP^d\dashrightarrow \Gr(d-r,r)\simeq\Gr(2r-d-1,r)$$
which associates to a general binary form $f$ of degree $d$ the projective
$(2r-d-1)$-dimensional 
 subspace 
 $\Pi_f=
 \PP((f^\perp)_r)\subset\PP(D_r)$ obtained from the degree $r$ component of the apolar ideal. 
In coordinates the map $\varPsi_{d,r}$ is defined by the maximal minors of the matrix 
$A^{d,r}$.
We denote by $Z_{d,r}=\overline{\varPsi_{d,r}(\PP^d)}\subset\Gr(2r-d-1,r)$ the closure of the image of $\varPsi_{d,r}$.

\section{Real rank boundaries of degree 7 binary forms}\label{secCaso7}
In this section, we prove the following:
\begin{theorem}\label{teorema}
  The real rank boundaries for degree $7$ binary real forms are the following hypersurfaces:
  \begin{align*}
    \partial_{\textrm{alg}}(\mathcal{R}_{7,4}) &= (\Delta_{(3,2,2)})^\vee; \\
  \partial_{\textrm{alg}}(\mathcal{R}_{7,5}) &=(\Delta_{(3,2,2)})^\vee\cup(\Delta_{(4,3)})^\vee\cup (\Delta_{(5,2)})^\vee ;\\
   \partial_{\textrm{alg}}(\mathcal{R}_{7,6}) &=(\Delta_{(4,3)})^\vee\cup (\Delta_{(5,2)})^\vee\cup(\Delta_{(7)})^\vee ; \\
  \partial_{\textrm{alg}}(\mathcal{R}_{7,7}) &=(\Delta_{(7)})^\vee .
  \end{align*}
\end{theorem}
\begin{remark}
From Proposition~\ref{propLSJoin} and formula \eqref{luke} we obtain:
\begin{itemize}
\item $(\Delta_{(3,2,2)})^\vee = \mathrm{Join}(\Delta_{(6,1)},\Delta_{(7)},\Delta_{(7)})$
is a hypersurface of  degree $24$;
\item $(\Delta_{(4,3)})^\vee = \mathrm{Join}(\Delta_{(5,1,1)},\Delta_{(6,1)})$
is a hypersurface of  degree $36$;
\item $(\Delta_{(5,2)})^\vee = \mathrm{Join}(\Delta_{(4,1,1,1)},\Delta_{(7)})$
is a hypersurface of  degree $24$;
\item $(\Delta_{(7)})^\vee = \Delta_{(2,1,1,1,1,1)}$ 
is a hypersurface of  degree $12$.
\end{itemize}
\end{remark}
\begin{proof} We divide the proof in several steps.

\medskip 
\noindent
\emph{The boundary $\partial_{\textrm{alg}}(\mathcal{R}_{7,7})$ between ranks 7 and $\leq 6$.}
\medskip
  
From \cite{ComonOttaviani} and \cite{CausaRe}, it is know that
the real rank boundary $\partial_{\textrm{alg}}(\mathcal{R}_{7,7})$ is
the discriminant hypersurface $\Delta_{(2,1^5)}$ in $\PP^7$.
Note that 
$$\Delta_{(2,1^5)}=(\Delta_{(7)})^\vee=
\overline{\varPsi_{7,6}^{-1}(CH_0(\Delta_{(6)}))}$$
where
$\varPsi_{7,6}: \PP^7\dashrightarrow Z_{7,6}\subset\Gr(4,6)\subset\PP^{20}$.

\medskip 
\noindent
\emph{The boundary \texorpdfstring{$\partial_{\textrm{alg}}(\mathcal{R}_{7,4})$}{D_alg(R_(7,4))} 
between ranks 4 and \texorpdfstring{$\geq5$}{>=5}.}
\medskip

    For the 
    reader's convenience, we sketch briefly the proof given in \cite{LeeSturmfels} of the fact that 
    $\partial_{\textrm{alg}}(\mathcal{R}_{7,4})= (\Delta_{(3,2,2)})^\vee$.
Consider a binary form $f$ of degree 7 with apolar ideal $f^\perp=(g_4,g_5)$, where $\deg(g_i)=i$.
By Lemma~\ref{apolarityLemma}, we have that $f\in \mathcal{R}_{7,4}$ if and only if $g_4$ is real-rooted.
When $f$ moves toward $\mathcal{R}_{7,5}\cup\mathcal{R}_{7,6}\cup\mathcal{R}_{7,7}$ and passes through the boundary
$\partial_{\textrm{alg}}(\mathcal{R}_{7,4})$, then (at least) two roots of $g_4$ 
must collapse and become a double root. 
Hence at the transition point the generator $g_4$ belongs to the discriminant locus
$\Delta_{(2,1,1)}\subset \PP(D_4)$, and 
by Proposition \ref{proposition-LS}, we get
$\partial_{\textrm{alg}}(\mathcal{R}_{7,4})\subseteq (\Delta_{(3,2,2)})^\vee$.
Now since $\partial_{\textrm{alg}}(\mathcal{R}_{7,4})\neq\emptyset$ because
$4$ and $5$ are typical ranks, and
  $\Delta_{(3,2,2)})^\vee$ is irreducible, it follows that
  $\partial_{\textrm{alg}}(\mathcal{R}_{7,4})= (\Delta_{(3,2,2)})^\vee$.

\medskip
\noindent
\emph{The boundary $\partial_{\textrm{alg}}(\mathcal{R}_{7,5})$ between ranks 5 and $\neq 5$.}
\medskip

We describe now
the boundary between 
$\mathcal{R}_{7,5}$ and $\mathcal{R}_{7,6}\cup\mathcal{R}_{7,7}$.
Let $f_{\varepsilon}$ be a continuous family of forms crossing the boundary 
$\partial_{\textrm{alg}}(\mathcal{R}_{7,5})$ 
at the point $f_0=f$, going 
from $\mathcal{R}_{7,5}$ to $\mathcal{R}_{7,6}\cup\mathcal{R}_{7,7}$. 
Namely, we assume that $f_{-\varepsilon}\in\mathcal{R}_{7,5}$ and 
$f_\varepsilon\in\mathcal{R}_{7,6}\cup \mathcal{R}_{7,7}$ for any small $\varepsilon$ with $\varepsilon>0$.
We can assume that for any $\varepsilon$ the form 
$f_\varepsilon$ is generated in generic degree, 
since the locus of non generated in generic degree forms has codimension $2$. 
In particular, we can assume $f^\perp =(g_4, g_5)$, where $\deg(g_4)=4$ and $\deg(g_5)=5$.

Let $\Gr(2,5)$ be the Grassmannian of planes in $\PP(D_5)$, and
consider the apolar map
\begin{equation}\label{apolar7}
\varPsi_{7,5}: \PP^7\dashrightarrow Z_{7,5}\subset\Gr(2,5)\subset\PP^{19},
\end{equation}
which is a cubic birational map onto a subvariety $Z_{7,5}\subset\PP^{19}$ of degree $84$
and cut out by $42$ quadric hypersurfaces.
The map $\varPsi_{7,5}$ sends the family $f_\varepsilon$ into 
a continuous family of apolar planes $\Pi_{\varepsilon}$; in particular,
$\varPsi_{7,5}(f)=\Pi_0=\langle ug_4,vg_4,g_5\rangle$ is the apolar plane of $f$.
From Lemma~\ref{apolarityLemma}, we obtain that the 
plane
$\Pi_{\varepsilon}$, with $\varepsilon<0$, 
contains a real-rooted form 
$h_\varepsilon=l_1(\varepsilon)l_2(\varepsilon)l_3(\varepsilon)l_4(\varepsilon)l_5(\varepsilon)$ 
(where $l_i(\varepsilon)\in D_1$), while 
$\Pi_{\varepsilon}$, with $\varepsilon>0$, does not contain any real-rooted form.
The set of real-rooted forms
is a full-dimensional connected  semi-algebraic subset of $\PP^5$, 
and the Zariski closure of its topological boundary is the discriminant hypersurface $\Delta=\Delta_{(2,1,1,1)}$.
Thus the limit 
$h_0=\lim_{\varepsilon\rightarrow 0^{-}} h_{\varepsilon} = l_1^2l_2l_3l_4$ must belong to 
$\Delta $.
We now 
 analyze, taking into account also Example~\ref{esempioSING},
 the possible positions of $\Pi_0$ with respect to $\Delta$:
\begin{enumerate}
\item\label{caso1di7}
The point $h_0$ is smooth 
and the tangent space
$T_{h_0}(\Delta) = \langle l_1 u^i v^j : i + j = 4 \rangle$ contains
$\Pi_0$.
This implies that $\Pi_0\in CH_2(\Delta)$.
\item\label{caso2}
  The point $h_0$ is smooth in a component of $\Delta_{(3,1,1)}\cup \Delta_{(2,2,1)}$. 
  We have the following subcases:
  \begin{enumerate}
  \item\label{caso2a}
    $h_0=l_1^3l_2l_3\in\Delta_{(3,1,1)}$ and
    $T_{h_0}(\Delta_{(3,1,1)}) = \langle l_1^2 u^i v^j : i + j = 3 \rangle$
    intersects $\Pi_0$ in a line $L$ through $h_0$.
    This implies that $\Pi_0\in CH_1(\Delta_{(3,1,1)})$.
    \item\label{caso2b}
        $h_0=l_1^2l_2^2l_3\in\Delta_{(2,2,1)}$ and
    $T_{h_0}(\Delta_{(2,2,1)}) = \langle l_1 l_2 u^i v^j : i + j = 3 \rangle$
      intersects $\Pi_0$ in a line $L$ through $h_0$.
      This implies that $\Pi_0\in CH_1(\Delta_{(2,2,1)})$.
  \end{enumerate}
\item\label{caso3}
  The point $h_0$ belongs to a component 
  of $\Delta_{(4,1)}\cup \Delta_{(3,2)}$, hence $\Pi_0\in CH_0(\Delta_{(4,1)})\cup CH_0(\Delta_{(3,2)})$. 
  \end{enumerate}
  
\noindent
Case~\eqref{caso1di7}.  Clearly this case cannot occur.
Indeed $g_4$ and $g_5$ would have $l_1$ as common divisor, 
and this is against our assumptions.

\medskip

\noindent
Case~\eqref{caso2}.
We show that case~\eqref{caso2a} cannot occur. 
With the same argument, 
one sees that neither case \eqref{caso2b} occurs.

If $h_0\not\in \langle ug_4,vg_4\rangle$, 
then we can take as degree 5 generator of the apolar ideal $g_5=h_0$.
Now, every point of $L$ is a form divisible by $l_1^2$
and we have that $L\cap \langle ug_4,vg_4 \rangle\neq\emptyset$. 
This implies that $l_1$ is a common divisor of $g_4$ and of $g_5$, which is impossible.
It follows that $h_0\in \langle ug_4,vg_4\rangle$ and in particular $l_1^2$ divides $g_4$.
More precisely, this implies that $g_4$ is of the form $l_1^2l_2l_3$, or $l_1^3l_2$, or $l_1^3l_3$. 
In any cases it is obvious that $f$ is limit of generic forms of degree $4$.
This implies that $f$ is a singular point of the hypersurface 
$\partial_{\textrm{alg}}(\mathcal{R}_{7,5})$. 
Hence $f$ does not vary in a codimension $1$ locus of $\PP^7$, and 
$\overline{\varPsi_{7,5}^{-1}(CH_1(\Delta_{3,1,1}))}$ 
cannot be a component of the boundary $\partial_{\textrm{alg}}(\mathcal{R}_{7,5})$.

\medskip
\noindent 
Case~\eqref{caso3}.
%
We show now that both components corresponding to this case 
are in the boundary.
Indeed it is enough to find an example of a binary form which lies exclusively on each component and is limit of a sequence of general forms of rank $5$ and a sequence of general forms of rank $6$.
This is done
in Example \ref{esempio7} below.
Recall that by Proposition \ref{proposition-LS}, we have 
$\overline{\varPsi_{7,5}^{-1}(CH_0(\Delta_{(4,1)}))}=(\Delta_{(5,2)})^\vee$
and
$\overline{\varPsi_{7,5}^{-1}(CH_0(\Delta_{(3,2)}))}=(\Delta_{(4,3)})^\vee$.
Hence we have proved that
$\overline{\partial_{\textrm{alg}}(\mathcal{R}_{7,5})\setminus \partial_{\textrm{alg}}(\mathcal{R}_{7,4})}
=(\Delta_{(4,3)})^\vee\cup (\Delta_{(5,2)})^\vee$.

\medskip
\noindent
\emph{The boundary $\partial_{\textrm{alg}}(\mathcal{R}_{7,6})$ between ranks 6 and $\neq 6$.}
\medskip

At this point we know that 
\begin{align*}
\overline{\partial_{\textrm{alg}}(\mathcal{R}_{7,6})\setminus \partial_{\textrm{alg}}(\mathcal{R}_{7,4})} &= 
\left(\overline{\partial_{\textrm{alg}}(\mathcal{R}_{7,5})\setminus \partial_{\textrm{alg}}(\mathcal{R}_{7,4})}\right) 
\cup \partial_{\textrm{alg}}(\mathcal{R}_{7,7})
\\ &=(\Delta_{(4,3)})^\vee\cup (\Delta_{(5,2)})^\vee\cup(\Delta_{(7)})^\vee.
\end{align*}
So, we only need to show that 
the boundary between $\mathcal{R}_{7,4}$ and $\mathcal{R}_{7,6}$ is not of codimension $1$ in $\mathbb{P}^7$.
Let $f_\varepsilon$ be a continuous family of forms such that $f_{-\varepsilon}\in\mathcal{R}_{7,4}$ and 
$f_\varepsilon\in\mathcal{R}_{7,6}$ for any small $\varepsilon$ with $\varepsilon>0$.
The corresponding apolar plane $\Pi_\varepsilon = \varPsi_{7,5}(f_\varepsilon)$
does not contain any real-rooted form for any $\varepsilon>0$. 
On the other hand, from Corollary~\ref{rangoreale}, we deduce that $\Pi_\varepsilon$
must contain a special line $L_\varepsilon$ 
which is generically contained
 in the locus of real-rooted forms for any  $\varepsilon<0$.
Now the limit $L_0=\lim_{\varepsilon\rightarrow 0^{-}} L_{\varepsilon}$
is a special line contained in the intersection of the plane $\Pi_{0}$ and of
the discriminant $\Delta = \Delta_{(2,1,1,1)}$.
By the previous analysis we deduce that the line $L_0=\langle ug_4,vg_4\rangle$ must be contained 
in $\Delta_{(4,1)}\cup\Delta_{(3,2)}$. This implies that $g_4\in\Delta_{(4)}$, and this 
forces $f_0$ to move in some locus of codimension $\ge2$  in $\mathbb{P}^7$, 
which cannot be a component of the boundary. 
\end{proof}
\begin{example}\label{esempio7}
Given
$$g_4=(u^2+v^2)(u^2- v^2),\quad g_5(\varepsilon)=(u^2+\varepsilon v^2)uv(\varepsilon u+v),$$
the degree $7$ form $f_{\varepsilon}$ associated to the apolar ideal $(g_4,g_5(\varepsilon))$ is:
\begin{multline*}
{\varepsilon}^{2} x^{7}+ 7 ({\varepsilon}^{2}+{\varepsilon}+1) x^{6} y -21 {\varepsilon} ({\varepsilon}^{2}+ {\varepsilon}+1) x^{5} y^{2} -35 {\varepsilon} x^{4} y^{3} +35 {\varepsilon}^{2} x^{3} y^{4} \\+ 21 ({\varepsilon}^{2}+{\varepsilon}+1) x^{2} y^{5}-7 {\varepsilon} ( {\varepsilon}^{2}+ {\varepsilon}+1) x y^{6}-{\varepsilon} y^{7} .
\end{multline*}
We have $\rk(f_{\varepsilon})=6$ for any small $\varepsilon\ge0$
and $\rk(f_{-\varepsilon})=5$ for any small $\varepsilon>0$. 
Moreover $f_0 = x^{6} y+3 x^{2} y^{5}$ belongs to 
$(\Delta_{(4,3)})^\vee$
and it does not belong to 
$(\Delta_{(5,2)})^\vee \cup  (\Delta_{(3,2,2)})^\vee$.

On the other hand, taking
$$g_4=(u^2+v^2)(2u^2- v^2),\quad g_5(\varepsilon)=(\varepsilon u^2+v^2)uv(\varepsilon u+v),$$
we consider the associated form $f_\varepsilon$:
\begin{multline*}
{\varepsilon} ({\varepsilon}^{3}+{\varepsilon}^{2}-{\varepsilon} -3) x^{7}+14({\varepsilon}^{2}- {\varepsilon}-1) x^{6} y -42 {\varepsilon}({\varepsilon}^{2}-{\varepsilon}-1) x^{5} y^{2} \\ -70 ({\varepsilon}^{3}-2) x^{4} y^{3}+70 {\varepsilon}({\varepsilon}^{3}-2) x^{3} y^{4} -42 {\varepsilon}({\varepsilon}^{2}-2 {\varepsilon}+2) x^{2} y^{5} \\ +  14 {\varepsilon}^{2} ({\varepsilon}^{2}-2 {\varepsilon}+2) x y^{6}-2(3 {\varepsilon}^{3}-2 {\varepsilon}^{2}+2 {\varepsilon}-4) y^{7} .
\end{multline*}
Again we have $\rk(f_{\varepsilon})=6$ for any small $\varepsilon\ge0$
and $\rk(f_{-\varepsilon})=5$ for any small $\varepsilon>0$. 
Moreover $f_0 = 7 x^{6} y-70 x^{4} y^{3}-4 y^{7}$ belongs to 
$(\Delta_{(5,2)})^\vee$ and it does not belong to 
$(\Delta_{(4,3)})^\vee \cup  (\Delta_{(3,2,2)})^\vee$.
For computational details, see Section~\ref{sezione-conti}. 
\end{example}

\section{Real rank boundaries of degree 8 binary forms}\label{secCaso8}
In this section, we prove the following:
\begin{theorem}\label{teorema2}
  The real rank boundaries for degree $8$ binary real forms are the following hypersurfaces:
  \begin{align*}
    \partial_{\textrm{alg}}(\mathcal{R}_{8,5}) &= (\Delta_{(3,3,2)})^\vee \cup (\Delta_{(4,2,2)})^\vee; \\
  \partial_{\textrm{alg}}(\mathcal{R}_{8,6}) &= (\Delta_{(3,3,2)})^\vee \cup (\Delta_{(4,2,2)})^\vee
  \cup(\Delta_{(4,4)})^\vee\cup (\Delta_{(5,3)})^\vee\cup (\Delta_{(6,2)})^\vee; \\
    \partial_{\textrm{alg}}(\mathcal{R}_{8,7}) &=
    (\Delta_{(4,4)})^\vee\cup (\Delta_{(5,3)})^\vee\cup (\Delta_{(6,2)})^\vee
    \cup (\Delta_{(8)})^\vee; \\
  \partial_{\textrm{alg}}(\mathcal{R}_{8,8}) &=(\Delta_{(8)})^\vee.
  \end{align*}
\end{theorem}
\begin{remark}
From Proposition~\ref{propLSJoin} and formula \eqref{luke} we obtain:
\begin{itemize}
\item $(\Delta_{(3,3,2)})^\vee = \mathrm{Join}(\Delta_{(7,1)},\Delta_{(7,1)},\Delta_{(8)})$
is a hypersurface of  degree $48$;
\item $(\Delta_{(4,2,2)})^\vee = \mathrm{Join}(\Delta_{(6,1,1)},\Delta_{(8)},\Delta_{(8)})$
is a hypersurface of  degree $36$;
\item $(\Delta_{(4,4)})^\vee = \mathrm{Join}(\Delta_{(6,1,1)},\Delta_{(6,1,1)})$
is a hypersurface of  degree $27$;
\item $(\Delta_{(5,3)})^\vee = \mathrm{Join}(\Delta_{(5,1,1,1)},\Delta_{(7,1)})$
is a hypersurface of  degree $48$;
\item $(\Delta_{(6,2)})^\vee = \mathrm{Join}(\Delta_{(4,1,1,1,1)},\Delta_{(8)})$
is a hypersurface of  degree $30$;
\item $(\Delta_{(8)})^\vee = \Delta_{(2,1,1,1,1,1,1)}$
is a hypersurface of  degree $14$.
\end{itemize}
\end{remark}
\begin{proof}
From \cite{ComonOttaviani} and \cite{CausaRe},
we have $\partial_{\textrm{alg}}(\mathcal{R}_{8,8})=(\Delta_{(8)})^\vee$. 
On the other hand, from \cite{LeeSturmfels}, 
we have $\partial_{\textrm{alg}}(\mathcal{R}_{8,5})= (\Delta_{(3,3,2)})^\vee \cup (\Delta_{(4,2,2)})^\vee$. 

We study now the boundary 
between ranks $6$ and $\ge6$.
Let $f_{\varepsilon}$ be a continuous family of forms crossing the boundary 
$\partial_{\textrm{alg}}(\mathcal{R}_{8,6})$ at the point $f_0=f$, going 
from $\mathcal{R}_{8,6}$ to $\mathcal{R}_{8,7}\cup\mathcal{R}_{8,8}$.
Namely, we assume that $f_{-\varepsilon}\in\mathcal{R}_{8,6}$ and 
$f_\varepsilon\in\mathcal{R}_{8,7}\cup \mathcal{R}_{8,8}$ for any small $\varepsilon$ with $\varepsilon>0$.
We can also assume that for any $\varepsilon\neq0$ the form $f_\varepsilon$ is generated in generic degree, \emph{i.e.}
the apolar ideal $f_{\varepsilon}^{\perp}$ is generated by two quintic forms
$g_\varepsilon$ and $g'_\varepsilon$.
Moreover, since $f$ 
moves in a codimension $1$ locus,
we may assume that 
$f^{\perp}$ is 
generated by two forms
$g_0$ and $g'_0$ either with $\deg(g_0) = \deg(g_0') = 5$,
or with $\deg(g_0) = 4$, $\deg(g_0')=6$, and moreover $g_0\not\in\Delta_{(2,1,1)}$.
Note that in the former case we have 
$\PP((f^{\perp})_6) = \langle ug_0,vg_0, ug'_0,vg'_0\rangle$,
while in the latter case we have 
$\PP((f^{\perp})_6) =  \langle u^2 g_0,uv g_0,v^2 g_0,g_0'\rangle$.

Consider the apolar map 
\begin{equation}\label{apolardegree8}
\varPsi_{8,6}: \PP^8\dashrightarrow Z_{8,6}\subset\Gr(3,6)\subset\PP^{34},
\end{equation}
which is a cubic birational map onto a subvariety $Z_{8,6}\subset\PP^{34}$ of degree $686$
and cut out by $186$ quadric hypersurfaces.
The map $\varPsi_{8,6}$ sends the family $f_\varepsilon$ into 
the continuous family of the $3$-dimensional linear spaces 
$\Pi_\varepsilon =\PP( (f_{\varepsilon}^{\perp})_6)$. 
From Lemma~\ref{apolarityLemma}, we obtain that 
$\Pi_{\varepsilon}$, with $\varepsilon<0$, 
contains a real-rooted form 
$h_\varepsilon=\prod_{i=1}^6 l_i(\varepsilon)$ (where $l_i\in D_1$), while 
$\Pi_{\varepsilon}$, with $\varepsilon>0$, does not contain any real-rooted form.
Thus the limit 
$h_0=\lim_{\varepsilon\rightarrow 0^{-}} h_{\varepsilon}$ must belong to the 
discriminant hypersurface $\Delta = \Delta_{(2,1,1,1,1)}$.
We now analyze, recalling Example~\ref{esempioSING}, the possible positions of $\Pi_0$ with respect to $\Delta$: 
\begin{enumerate}
\item\label{caso8-1}
The point $h_0$ is smooth  and the tangent space 
$T_{h_0}(\Delta) = \langle l_1 u^i v^j : i + j = 5 \rangle$ 
contains
$\Pi_0$.
This implies that $\Pi_0\in CH_3(\Delta)$.
\item\label{caso8-2}
  The point $h_0$ is smooth in a component 
  of $\Delta_{(3,1,1,1)}\cup \Delta_{(2,2,1,1)}$. 
We have the following subcases:
  \begin{enumerate}
  \item\label{casi8CH2-a}
    $h_0=l_1^3l_2l_3l_4\in\Delta_{(3,1,1,1)}$ and
     $T_{h_0}(\Delta_{(3,1,1,1)}) = \langle l_1^2 u^i v^j : i + j = 4 \rangle$ 
    intersects $\Pi_0$ in a plane $P$ through $h_0$.
    This implies that $\Pi_0\in CH_2(\Delta_{(3,1,1,1)})$.
    \item\label{casi8CH2-b}
        $h_0=l_1^2l_2^2l_3l_4\in\Delta_{(2,2,1,1)}$ and
        $T_{h_0}(\Delta_{(2,2,1,1)}) = \langle l_1 l_2 u^i v^j : i + j = 4 \rangle$ 
      intersects $\Pi_0$ in a plane $P$ through $h_0$.
      This implies that $\Pi_0\in CH_2(\Delta_{(2,2,1,1)})$.
  \end{enumerate}
  \item\label{casi8CH1}
  The point $h_0$ is smooth in a component of 
  $\Delta_{(3,2,1)} \cup \Delta_{(4,1,1)} \cup \Delta_{(2,2,2)}$.
  We have the following subcases:
  \begin{enumerate}
  \item\label{casi8CH1-a}
    $h_0=l_1^3l_2^2l_3\in\Delta_{(3,2,1)}$ and
    $T_{h_0}(\Delta_{(3,2,1)}) = \langle l_1^2 l_2 u^i v^j : i + j = 3 \rangle$
    intersects $\Pi_0$ in a line $L$ through $h_0$.
    This implies that $\Pi_0\in CH_1(\Delta_{(3,2,1)})$.
  \item \label{casi8CH1-b}
    $h_0=l_1^4l_2l_3\in\Delta_{(4,1,1)}$ and
    $T_{h_0}(\Delta_{(4,1,1)}) = \langle l_1^3 u^i v^j : i+j = 3 \rangle$
    intersects $\Pi_0$ in a line $L$ through $h_0$.
    This implies that $\Pi_0\in CH_1(\Delta_{(4,1,1)})$.
  \item \label{casi8CH1-c}
    $h_0=l_1^2l_2^2l_3^2\in\Delta_{(2,2,2)}$ and
    $T_{h_0}(\Delta_{(2,2,2)}) = \langle l_1 l_2 l_3 u^i v^j: i + j = 3 \rangle$
    intersects $\Pi_0$ in a line $L$ through $h_0$.
    This implies that $\Pi_0\in CH_1(\Delta_{(2,2,2)})$.
  \end{enumerate}
\item\label{componentivere}
  The point $h_0$ belongs to a component of 
  $\Delta_{(3,3)}\cup \Delta_{(4,2)} \cup \Delta_{(5,1)}$, hence 
  $\Pi_0\in CH_0(\Delta_{(3,3)})\cup CH_0(\Delta_{(4,2)}) \cup CH_0(\Delta_{(5,1)})$. 
  \end{enumerate}
  In the following, we show that only the last case occurs.
  
\medskip
\noindent 
Case \eqref{caso8-1}.
This case cannot occur. Indeed from the fact that 
$\Pi_0\subset T_{h_0}(\Delta)$,
we would conclude that $l_1$ is a common divisor of $g_0$ and $g_0'$. 

\medskip
\noindent 
Case \eqref{caso8-2}. 
Consider first the case when $f$ is not generated in generic degree, 
so that we have 
$\Pi_0 = \langle u^2 g_0,uv g_0,v^2 g_0,g_0'\rangle$.
  If $h_0\not\in \langle u^2 g_0,uv g_0,v^2 g_0\rangle$, then we can take $g_0'=h_0$ and, 
since there are at least two points in $P\cap \langle u^2 g_0,uv g_0,v^2 g_0\rangle$, 
we deduce that $g_0$ and $g_0'$ have a common divisor, which is a contradiction.
  Assume therefore that $h_0\in \langle u^2 g_0,uv g_0,v^2 g_0\rangle$. 
  Since 
  $g_0\not\in\Delta_{(2,1,1)}$, 
  the only possibility is that $g_0=l_1l_2l_3l_4$. 
  This implies that $\rk(f)=4$, and it is easy to see that $f$ is limit of a general 
  sequence of form of rank $5$. This would implies that $f$ is not a general point of the 
  boundary between forms of rank $6$ and rank $\geq7$.
  Hence we can assume 
  $\deg(g_0)=\deg(g_0')=5$, and consider the following subcases.

  \medskip
\noindent 
Case \eqref{casi8CH2-a}. The plane $P$ meets the special lines 
$\langle ug_0,vg_0\rangle$ and
$\langle ug'_0,vg'_0\rangle$
at points
$p_0$ and $p'_0$ respectively. 
Therefore, $p_0$ and $p'_0$ are forms
 divisible by $l_1^2$, and then $l_1$ divides 
both $g_0$ and $g'_0$, which is a contradiction.

\medskip
\noindent
Case \eqref{casi8CH2-b}. 
Let us consider the surface
$$\mathcal{Q}=\bigcup_{\begin{subarray}{c} m\in D_1, \\ g\in \langle g_0,g_0'\rangle\end{subarray}} m g \subset \Pi_0\simeq \mathbb{P}^3$$
swept out by all the special apolar lines of $f$.
Using that $f$ is generated in generic degrees, 
one sees that $\mathcal{Q}$ is  a smooth quadric surface, which we will call the \emph{apolar quadric} of $f$.
The intersection $P\cap \mathcal{Q}$ is a (possible reducible) plane conic, which 
in particular contains three noncollinear points: 
$p_0 = m\,g$, $p'_0 = m'\,g'$ and $p''_0 = m''\,g''$.
We can assume $g=g_0$, $g'=g_0'$, 
and since every point of $P$ is a form divisible by $l_1\,l_2$, 
we conclude that $g_0$ and $g_0'$ have a common factor, which is a contradiction.

\medskip
\noindent 
Case \eqref{casi8CH1}.
As above we consider first the  case 
when  $\deg(g_0)=4$, $\deg(g_0')=6$ and $\Pi_0 = \langle u^2 g_0,uv g_0,v^2 g_0,g_0'\rangle$.
  If $h_0\not\in \langle u^2 g_0,uv g_0,v^2 g_0\rangle$, then we can take $g_0'=h_0$.
  Moreover, since $L\cap \langle u^2 g_0,uv g_0,v^2 g_0\rangle\neq\emptyset$,
we deduce that $g_0$ and $g_0'$ have a common divisor, which is a contradiction.
Thus we have that $h_0\in \langle u^2 g_0,uv g_0,v^2 g_0\rangle$. This implies 
that $g_0\in\Delta_{(2,1,1)}$,
which contradicts our assumption.
Hence we can assume 
$\deg(g_0)=\deg(g_0')=5$, and consider the following subcases.

\medskip
\noindent
Case \eqref{casi8CH1-a}.
Let $\mathcal{Q}$ 
be again the apolar quadric of $f$.
We have two cases: either the line $L$ meets $\mathcal{Q}$ in two distinct points $m g$ and $m' g'$, 
or there exists a point $m g\in L\cap \mathcal{Q}$ such that $L$ is contained in the tangent plane
$T_{mg}\mathcal{Q}$.

In the former case, since $l_1^2 l_2$ divides $m g$ and $m' g'$, 
we deduce that
$l_1$ divides $g$ and $g'$. This is a contradiction,
unless we have $g = g'$ and hence $L$ is the special line $\langle g u, g v\rangle$.
Now,
since $h_0\in L$, we obtain that $g\in \Delta_{(2,2,1)}\cup\Delta_{(3,1,1)}$, 
 and thus 
 $f\in (\Delta_{(3,3,2)})^{\vee} \cup (\Delta_{(4,2,2)})^{\vee} = \partial_{\textrm{alg}}(\mathcal{R}_{8,5})$
 is also limit of generic forms of rank $5$.
This implies that $f$ belongs to the singular locus of the hypersurface 
$\partial_{\textrm{alg}}(\mathcal{R}_{7,6})$ and then
$\overline{\varPsi_{8,6}^{-1}(CH_1(\Delta_{(3,2,1)}))}$ 
cannot be a component of the boundary $\partial_{\textrm{alg}}(\mathcal{R}_{8,6})$.

In the latter case, 
we may assume $g = g_0$ and $T_{mg}\mathcal{Q} = \langle m g_0, m g_0', m' g_0 \rangle$, for some $m'\in D_1$.
We have $h_0 = l_1^3l_2^2l_3 = \alpha m g_0 + \beta m g_0' + \gamma m' g_0$, for some scalars $\alpha,\beta,\gamma$, 
and we know that $l_1^2 l_2$ divides $m g_0$. Since $\gcd(g_0,g_0') = 1$, this implies $\beta = 0$,
and then $l_1^3 l_2^2 l_3 = (\alpha m + \gamma m') g_0$. As above, from this it follows 
that $g_0\in \Delta_{(2,2,1)}\cup\Delta_{(3,1,1)}$ and thus 
$f$ does not vary in a codimension $1$ locus of $\PP^8$,

\medskip
\noindent 
Cases \eqref{casi8CH1-b} and 
\eqref{casi8CH1-c}.
Arguing as above, we deduce that $f$ must belong to 
$ (\Delta_{(4,2,2)})^{\vee}$ and
$(\Delta_{(3,3,2)})^{\vee}$,
respectively, and furthermore we must have that 
$f$ is a singular point of the hypersurface $\partial_{\textrm{alg}}(\mathcal{R}_{8,6})$.
This implies 
that $\overline{\varPsi_{8,6}^{-1}(CH_1(\Delta_{(4,1,1)}))}$ and 
$\overline{\varPsi_{8,6}^{-1}(CH_1(\Delta_{(2,2,2)}))}$ 
are not components of the boundary $\partial_{\textrm{alg}}(\mathcal{R}_{8,6})$.

\medskip
\noindent
Case~\eqref{componentivere}. 
We show in Example \ref{esempio8} below that each of the three components corresponding 
  to this case  are in the boundary.
This proves that
$\overline{\partial_{\textrm{alg}}(\mathcal{R}_{8,6})\setminus \partial_{\textrm{alg}}(\mathcal{R}_{8,5})}= (\Delta_{(4,4)})^\vee\cup (\Delta_{(5,3)})^\vee\cup (\Delta_{(6,2)})^\vee$.

\medskip 
Finally we need to prove that there are no components of the boundary 
between $\mathcal{R}_5$ and $\mathcal{R}_7$. 
This can be done with the same argument used at the end of the proof of Theorem \ref{teorema}.
\end{proof}

\begin{example}\label{esempio8}
 Given
 \[ g_0 = u^{4} v-u^{2} v^{3}-2 v^{5}, \quad g_0' = -u^{5}+2 u^{3} v^{2}+2 u v^{4},\]
we have $u g_0 + v g_0' = u^3 v^3$ and
the degree $8$ form $f_0$ associated to the apolar ideal $(g_0,g_0')$ is 
\begin{equation}\label{primaFdegree8}
f_0 = 8 x^{8}+112 x^{6} y^{2}+56 x^{2} y^{6}-y^{8}. 
\end{equation}
With the help of a computer, one can easily check (see Section~\ref{sezione-conti}) that
$\rk(f_0)=7$ and $f_0\in (\Delta_{(4,4)})^\vee\setminus\left((\Delta_{(5,3)})^\vee\cup (\Delta_{(6,2)})^\vee\right)$. 
Moreover, 
we can construct near $f_0$ generic degree $8$ forms $f_{\pm\varepsilon}$ having real ranks $6$ and $7$. 

Analogously, given
\[g_0 = u^{4} v-u^{2} v^{3}-2 v^{5}, \quad g_0' = -u^{5}+u^{4} v+u^{3} v^{2}+2 u v^{4},\]
we have $u g_0 + v g_0' = u^4 v^2$ 
and the associated degree $8$ form is
\begin{equation}
f_0 = x^{8}+8 x^{7} y+28 x^{3} y^{5}-2 x y^{7} . 
\end{equation}
 One verifies 
 that $\rk(f_0)=7$ and $f_0\in (\Delta_{(5,3)})^\vee\setminus\left((\Delta_{(4,4)})^\vee\cup (\Delta_{(6,2)})^\vee\right)$. 

 Finally, given 
\[  g_0 = u^{5}+u^{4} v+3 u^{3} v^{2}+3 u^{2} v^{3}+2 u v^{4}+2 v^{5}, \quad g_0' = -3 u^{3} v^{2}-2 u v^{4}, \]
 we have $u g_0 + (u+v) g_0' = u^5 (u+v)$ 
 and the associated degree $8$ form is 
 
\begin{multline}
 f_0 = 8 x^{8}-64 x^{7} y+224 x^{6} y^{2}-448 x^{5} y^{3}-840 x^{4} y^{4} \\ +672 x^{3} y^{5}+504 x^{2} y^{6}-144 x y^{7}-17 y^{8} .
\end{multline}
 One verifies 
 that $\rk(f_0)=7$ and $f_0\in (\Delta_{(6,2)})^\vee\setminus\left((\Delta_{(4,4)})^\vee\cup (\Delta_{(5,3)})^\vee\right)$. 
\end{example}

\section{Computations}\label{sezione-conti}\label{secComp}
%
%
We provide a package for \textsc{Macaulay2} \cite{macaulay2}, 
named 
\href{http://www2.macaulay2.com/Macaulay2/doc/Macaulay2-1.12/share/doc/Macaulay2/CoincidentRootLoci/html/index.html}{\texttt{CoincidentRootLoci}}
and 
included with the current stable version of \textsc{Macaulay2},
which implements methods useful to check the correctness of Examples~\ref{esempio7} and \ref{esempio8}.
This package depends on the packages \texttt{Cremona} and \texttt{Resultants}
(see \cite{packageCremona} and \cite{packageResultants}).
In the following, we illustrate briefly some of the methods available.
For technical details and examples, we refer to the documentation of the package, which
can be shown using the command \texttt{viewHelp}.

The method \texttt{realrank} computes the real rank of a binary 
form with rational coefficients. Indeed,
Lemma~\ref{apolarityLemma} reduces the problem of computing the real rank of a binary form to that 
of establishing whether certain semi-algebraic sets are nonempty.
The Tarski formulas defining these semi-algebraic 
sets can be obtained 
via the computation of kernels of appropriate 
catalecticant matrices. 
The problem of deciding the truth of a Tarski formula
can be handled by \textsc{Qepcad B} 
via a quantifier elimination by partial cylindrical algebraic decomposition (see \cite{qepcadBrown}).
The method  calls automatically \textsc{Qepcad B}
without requiring user intervention (provided it is installed on the system).
Below, we compute the real rank of the binary form \eqref{primaFdegree8} 
(the run time is about 30 seconds).
{\footnotesize
\begin{Verbatim}[commandchars=&!$]
Macaulay2, version 1.12
with packages: &colore!airforceblue$!ConwayPolynomials$, &colore!airforceblue$!Elimination$, &colore!airforceblue$!IntegralClosure$, &colore!airforceblue$!InverseSystems$, 
               &colore!airforceblue$!LLLBases$, &colore!airforceblue$!PrimaryDecomposition$, &colore!airforceblue$!ReesAlgebra$, &colore!airforceblue$!TangentCone$
&colore!darkorange$!i1 :$ &colore!airforceblue$!needsPackage$ "&colore!bleudefrance$!CoincidentRootLoci$"; 
&colore!darkorange$!i2 :$ R := &colore!darkspringgreen$!QQ$[x,y];
&colore!darkorange$!i3 :$ F = 8*x^8+112*x^6*y^2+56*x^2*y^6-y^8;
&colore!darkorange$!i4 :$ &colore!bleudefrance$!realrank$ F
o4 = 7
\end{Verbatim}
} \noindent 

The method \texttt{member} tests membership of a binary form  
in the dual variety of a coincident root locus (or in a coincident root locus).
It does not pass through the hard computation of the equations but uses Proposition~\ref{proposition-LS}.
Below, we verify that the binary form \eqref{primaFdegree8} lies in
$(\Delta_{(4,4)})^\vee$ but not in $(\Delta_{(5,3)})^\vee\cup (\Delta_{(6,2)})^\vee$ (the run time is less than one second).
{\footnotesize
\begin{Verbatim}[commandchars=&!$]
&colore!darkorange$!i5 :$ X = &colore!bleudefrance$!dual coincidentRootLocus$(4,4)
o5 = CRL(6,1,1) * CRL(6,1,1) (dual of CRL(4,4))
o5 : JoinOfCoincidentRootLoci
&colore!darkorange$!i6 :$ &colore!bleudefrance$!member$(F,X)
o6 = true
&colore!darkorange$!i7 :$ &colore!bleudefrance$!member$(F,&colore!bleudefrance$!dual coincidentRootLocus$(5,3)) &colore!darkorchid$!or$ &colore!bleudefrance$!member$(F,&colore!bleudefrance$!dual coincidentRootLocus$(6,2))
o7 = false
\end{Verbatim}
} \noindent 

The method 
\texttt{apolar} computes the apolar ideal of a binary form, while \texttt{recover},
as the name suggests, recovers the binary form from its apolar ideal. 
Basically, these two methods translate to problems of computing the 
image or the inverse image of a point via a (bi)rational map, and then the computation 
is performed using tools of the package \texttt{Cremona}.
For example, the following calculation involves the birational map \eqref{apolardegree8} 
(the run time is less than one second).
{\footnotesize
\begin{Verbatim}[commandchars=&!$]
&colore!darkorange$!i8 :$ F == &colore!bleudefrance$!recover apolar$ F
o8 = true
\end{Verbatim}
} \noindent 

For the convenience of the user, the method \texttt{realRankBoundary} implements 
Theorems~\ref{teorema} and \ref{teorema2}. For example, below we get immediately 
the degree of the first component of $\partial_{\textrm{alg}}(\mathcal{R}_{8,5})$.
{\footnotesize
\begin{Verbatim}[commandchars=&!$]
&colore!darkorange$!i9 :$ Y = &colore!airforceblue$!first$ &colore!bleudefrance$!realRankBoundary$(8,5)
o9 = CRL(7,1) * CRL(7,1) * CRL(8) (dual of CRL(3,3,2))
o9 : JoinOfCoincidentRootLoci
&colore!darkorange$!i10 :$ &colore!bleudefrance$!degree$ Y
o10 = 48
\end{Verbatim}
} \noindent 


\providecommand{\bysame}{\leavevmode\hbox to3em{\hrulefill}\thinspace}
\providecommand{\MR}{\relax\ifhmode\unskip\space\fi MR }
\providecommand{\MRhref}[2]{%
  \href{http://www.ams.org/mathscinet-getitem?mr=#1}{#2}
}
\providecommand{\href}[2]{#2}

\end{document}